\documentclass{amsart}
\usepackage[curve]{xypic}
\newbox\mybox
\def\overtag#1#2#3{\setbox\mybox\hbox{$#1$}\hbox to
  0pt{\vbox to 0pt{\vglue-#3\vglue-\ht\mybox\hbox to \wd\mybox
      {\hss$\ss#2$\hss}\vss}\hss}\box\mybox}
\def\undertag#1#2#3{\setbox\mybox\hbox{$#1$}\hbox to 0pt{\vbox to
    0pt{\vglue#3\vglue\ht\mybox\hbox to \wd\mybox
      {\hss$\ss#2$\hss}\vss}\hss}\box\mybox}
\def\lefttag#1#2#3{\hbox to 0pt{\vbox to 0pt{\vss\hbox to
      0pt{\hss$\ss#2$\hskip#3}\vss}}#1}
\def\righttag#1#2#3{\hbox to 0pt{\vbox to 0pt{\vss\hbox to
      0pt{\hskip#3$\ss#2$\hss}\vss}}#1}
\let\ss\scriptstyle

\def\Dot{\lower.2pc\hbox to 2pt{\hss$\bullet$\hss}}
\def\Circ{\lower.2pc\hbox to 2pt{\hss$\circ$\hss}}
\def\Vdots{\raise5pt\hbox{$\vdots$}}
\newcommand\lineto{\ar@{-}}
\newcommand\dashto{\ar@{--}}
\newcommand\dotto{\ar@{.}}
\newcommand\R{\mathbb R}
\newcommand\C{\mathbb C}
\newtheorem{theorem}{Theorem}[section]
\newtheorem*{question*}{Question}
\begin{document}
\title{Topology of Hypersurface Singularities}% (PRELIMINARY VERSION!)}
\author{Walter D. Neumann} \thanks{Supported under NSF grant no.\ 
  DMS-0083097. The hospitality of the Math. Research Institute in
  Bordeaux during the writing of this paper is gratefully
  acknowledged, as are comments by P. Cassou-Nogu\`es and A. Durfee on
  an early version. Durfee's article \cite{durfee99} in ``History of
  Topology'' is a recommended complement to this one.}
\address{Department of Mathematics\\Barnard College, Columbia
  University\\New York, NY 10027} \email{neumann@math.columbia.edu}
\keywords{plane curve singularity, link of a singularity, splice
  diagram}
\subjclass{}%14H20, 32S50, 57M25}
\begin{abstract}
  K\"ahler's paper \emph{\"Uber die Verzweigung einer algebraischen
    Funktion zweier Ver\"anderlichen in der Umgebung einer
    singul\"aren Stelle}'' offered a more perceptual view of the link
  of a complex plane curve singularity than that provided shortly
  before by Brauner. K\"ahler's innovation of using a ``square
  sphere'' became standard in the toolkit of later researchers on
  singularities. We describe his contribution and survey developments
  since then, including a brief discussion of the topology of isolated
  hypersurface singularities in higher dimension.
\end{abstract}
\maketitle 

\section{Topology of plane curve singularities}

The Riemann surface of an algebraic function on the plane represents a
complex curve (real dimension 2) as a covering of the Riemann sphere,
ramified over some finite collection of points. At the start of the
20th century, the study of complex surfaces (real dimension 4) was
rapidly developing, and they too were often studied as ``Riemann
surfaces,'' --- now of algebraic functions on the complex plane.  The
branching of such a ``Riemann surface'' is along a complex curve, and
the only difficult case in understanding the local topology of this
branching is at a singularity of the curve.  The problem therefore
arose, to understand the topology of a complex plane curve $C$ near a
singular point. 

The first discussion of this appears to be in
Heegaard's 1898 thesis \cite{heergaard} (see Epple \cite{epple,epple1}).
A small ball $B$ around the singular point will intersect the curve
$C$ in a set that is homeomorphic to the cone on $C\cap \partial B$.
This set $C\cap\partial B$, which is a link (disjoint union of
embedded circles) in the 3-sphere $\partial B$, therefore determines
the local topology completely.  It thus suffices to understand the
links that arise this way: \emph{links of plane curve singularities},
as they are now called. 

To understand the local branching of the ``Riemann surface'' one also
needs the fundamental group of the complement of the link in the
3-sphere.  The first comprehensive article on this topic is the 1928
paper \cite{brauner} by Karl Brauner, who writes that he learned the
problem from Wirtinger, who had spoken on it to the
Mathematikervereinigung in Meran in 1905 and subsequently held a
seminar on the topic in Vienna.  In his paper Brauner follows
Heegaard in using stereographic projection to move the link from
$S^3$ to $\R^3$. He then describes the topology of the link in terms
of repeated cabling, and gives an explicit presentation of the
fundamental group of the complement of the link.

Brauner's exposition is complicated by the stereographic projection,
and Erich K\"ahler revisits the question in the article
``\emph{\"Uber die Verzweigung einer algebraischen Funktion zweier
  Ver\"anderlichen in der Umgebung einer singul\"aren Stelle}''
\cite{kaehler}, pro\-mising and providing a more perceptual view than
Brauner's.  After choosing local coordinates centered at the singular
point, he replaces the round sphere $\partial B$ by a ``rectangular
sphere''
$$\{(x,y)\in \C^2: |x|=\epsilon, |y|\le \delta\text{~~or~~}
|x|\le\epsilon, |y|=\delta\}=\partial(D^2(\epsilon)\times
D^2(\delta)).$$
By choosing the coordinates suitably he also arranges
that the curve meets this rectangular sphere only in the portion
$|x|=\epsilon$. Since this portion is a solid torus that can be
identified with a standard solid torus in $\R^3$, this makes the
topology easier to visualize. Nowadays his technique is used
routinely, but it is reasonable to guess that the timespan between
Wirtinger's seminar and a general description of the topology was in
part due to the lack of an easy visualization technique.

To be more specific about K\"ahler's approach, suppose our curve is
given in local coordinates by an equation $f(x,y)=0$.  Newton had
already pointed out long before that one can give approximate
solutions to this equation, giving $y$ in terms of fractional powers
of $x$. Assuming the $y$-axis is not tangent to the curve at $(0,0)$,
Newton's successive approximations have the form
\begin{align*}
  y&=a_1x^{\frac{q_1}{p_1}}\\
  y&=x^{\frac{q_1}{p_1}}(a_1+a_2x^{\frac{q_2}{p_1p_2}})\\
  y&=x^{\frac{q_1}{p_1}}(a_1+x^{\frac{q_2}{p_1p_2}}(a_2+a_3x^{\frac{q_2}{p_1p_2p_3}}))\\
&\quad\vdots
\end{align*}
with $p_i$ and $q_i$ relatively prime positive integers.  There may be
several solutions of this type near $(x,y)=(0,0)$, corresponding to
different branches of the curve at $(0,0)$. The fact that the curve is
not tangent to the $y$-axis implies that $q_1\ge p_1$ for each branch,
and by choosing $\delta$ and $\epsilon$ suitably (they should be
small, and $\delta/\epsilon$ should exceed the absolute value of
the coefficient $a_1$ for each branch with $q_1=p_1$), one arranges
that the curve intersects only the solid torus $x=\epsilon$ of
K\"ahler's rectangular sphere.

It is now easy to see that the first Newton approximation gives a link
that is a $(p_1,q_1)$ torus knot, represented by a closed braid with
$p_1$ strands in this solid torus. The next approximation replaces
this by the $(p_2,q_2)$-cable\footnote{When talking of a $(p,q)$ cable
  on a knot, $q$ is only well defined after choosing a framing of the
  knot, that is, a choice of a parallel copy to call the $(1,0)$
  cable.  The framing we are using here is the ``naive'' framing of a
  cable knot, determined by choosing the parallel copy on the same
  torus that the cable knot naturally lies on. We return to the
  framing issue later.}  on this knot, represented by a
$p_1p_2$-strand braid, and so on. Thus each branch of the curve leads
to a component of the link that is an iterated cabling on a torus
knot. Such a link is called an \emph{iterated torus link}.

K\"ahler actually used the more familiar expression
of Newton's approximations as the partial sums of a fractional power
series solution
$$y=b_1x^\frac{r_1}{p_1}+b_2x^\frac{r_2}{p_1p_2}+
b_3x^\frac{r_3}{p_1p_2p_3}+\dots~$$
to $f(x,y)=0$, which had been
introduced by Puiseux in the mid-nineteenth century. The
pairs $(p_i,r_i)$ occurring in the exponents of the Puiseux series are
called \emph{Puiseux pairs}.  They of course determine and
are determined by the \emph{Newton pairs} $(p_i,q_i)$. The precise
inductive relationship is $q_1=r_1$, $q_i=r_i-p_ir_{i-1}$.

Not all Puiseux pairs are topologically significant\footnote{The
  terminology \emph{characteristic pair} is often used to single out
  the topologically significant Puiseux pairs, although this is with
  hindsight; the characteristic pairs were originally singled out for
  geometric reasons.}: a pair with $p_i=1$ does not contribute to the
topology of its link component, since it represents a $(1,q)$ cabling
for some $q$, which simply replaces a knot by a parallel copy of
itself. However this pair may, nevertheless, be topologically
significant, in that it \emph{can} contribute to the linking of
different link components with each other.  Thus care must be taken in
attempting to retain only the topologically significant data.
K\"ahler satisfies himself with describing typical cases that must be
considered in an iterative understanding of the topology and
fundamental group of any given example, but he gave no general
solution to this issue, writing: ``Es soll uns jedoch gen\"ugen, an
den vorstehenden bereits sehr allgemeinen Beispielen die
merkw\"urdigen Verzweigungsverh\"altnisse der Funktionen mehrerer
Variablen dargetan zu haben'' (It should suffice to have presented the
remarkable branching behavior of functions of several variables by
these already very general examples).

Although K\"ahler's presentation indeed provides the techniques to
deal with any particular example, it gives no explicit closed form
encapsulation of the topology. A question that was therefore addressed
by many later authors was:

\begin{question*}
  What invariant or collection of invariants completely determines the
  topology of a plane curve singularity?
\end{question*}

The implicit answer of K\"ahler's paper
is simply to retain relevant parts of the Puiseux expansions for each
branch, where ``relevant'' can be taken to mean: whenever two branches
have identical Puiseux expansions up to some point, include the final
term where they agree, and otherwise include only those terms that are
topologically relevant to a branch.  

A classical notion of equivalence of plane curve singularities, which,
with hindsight, is the same as topological equivalence, is based on
the the tree of infinitely near points, or, what is essentially
equivalent, the resolution diagram (see %Enriques-Chisini
\cite{enriques}). Various characterizations of this equivalence are
given in Zariski's investigation of equisingularity \cite{zariski1}. A
classical characterization, according to Reeve \cite{reeve}, is that
an equivalence class is determined by the sequence of characteristic
Puiseux pairs for each branch and the pairwise intersection numbers of
the branches (for a proof see \cite{zariski71} or \cite{le-jeune}).
Reeve shows that these intersection numbers are the linking numbers of
the corresponding components of the link of the singularity.  Thus:
\begin{theorem}
\label{th:1}
The link of a plane curve singularity is determined by 
%the topologies of the individual components, as coded by their 
the sequences of characteristic Puiseux pairs of the individual
components, and their pairwise linking numbers.
\end{theorem}

This presupposes agreement of classical and topological equivalence,
first proved for one branch in 1932 by Burau \cite{burau1} and Zariski
\cite{zariski0} independently, and then for two branches in 1934 by
Burau \cite{burau2}, who points out that the general case follows.
The connection with classical equivalence is not explicit in
\cite{burau2}, but was presumably understood. It is explicit in
Reeve's exposition.  Both Burau and Zariski recover the Puiseux data
for the link from its Alexander polynomial, and they use that the link
is the link of a plane curve singularity. In 1953 H. Schubert showed
that one can unravel the cabling numbers from the topology
of any cabled link \cite{schubert}.

Generalizing the work of Burau and Zariski, Evers
\cite{evers} and Yamamoto \cite{yamamoto} independently showed:
\begin{theorem}
  The multi-variable Alexander polynomial of the link of a plane curve
  singularity is a complete invariant for its topology.
\end{theorem}

R. Waldi, in his Regensburg dissertation \cite{waldi} showed:
\begin{theorem}
  The value semigroup of a plane curve singularity is a complete
  invariant for its topology.
\end{theorem}

Attractive as the above results are, they are not entirely
satisfactory as an encoding of the topology:
each does so %encode the topology 
in terms of a
redundant set of data from which other useful invariants are not
necessarily easy to compute.

In the 1970's there was a revolution in 3-manifold topology, brought
in part\footnote{The other part was Thurston's geometrization
  conjecture, which is, however, irrelevant to the 3-manifolds that
  arise in algebraic geometry.} by the JSJ decomposition theorem for
3-manifolds (foreseen by Waldhausen \cite{waldhausen} in a
little-noticed paper, and proved by Jaco-Shalen \cite{jaco-shalen} and
Johannson \cite{johannson}). In particular, this canonical
decomposition of any 3-manifold provides a general framework for (and
radical generalization of) Schubert's results for links mentioned
above, and hence a new view of the fact that classical and 
topological equivalence of plane curve singularities are the same.

In the early 1980's Eisenbud and the author used JSJ decomposition to
provide a new combinatorial encoding of
%the data relevant to 
the topology of a plane curve singularity: the
\emph{splice diagrams} of \cite{eisenbud-neumann}
(adapted\footnote{Actually, a case of convergent evolution.} from a
concept due to Siebenmann \cite{siebenmann}).  
On a superficial level, the splice diagram is
just another way of encoding the cabling information,
i.e., the Puiseux data.  The Puiseux pairs are replaced by new pairs
that have global topological meaning. For instance a single branch
with Puiseux pairs $(p_1,r_1), (p_2,r_2),\dots, (p_k,r_k)$ is encoded
by a splice diagram
$$
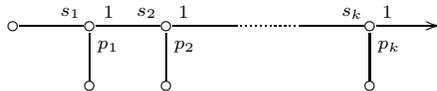
\xymatrix@R=18pt@C=24pt@M=0pt@W=0pt@H=0pt{
\circ\lineto[r]^(.75){s_1}&
\circ\lineto[d]^(.35){p_1}\lineto[r]^(.25){1}^(.75){s_2}&
\circ\lineto[d]^(.35){p_2}\lineto[r]^(.25){1}&\dotto[r]&\lineto[r]^(.75){s_k}&
\circ\lineto[d]^(.35){p_k}\ar[r]^(.25){1}&\\
&\circ&\circ&&&\circ}
$$
with $s_1=r_1$ and, for $i\ge1$, $s_{i+1}=r_{i+1}-r_ip_{i+1}+
p_ip_{i+1}s_i$.  The pairs $(p_i,s_i)$ describe the repeated cabling
in terms of the natural topological framings of  knots\footnote{This is
  the framing in which a parallel copy of a knot has zero linking
  number with the knot.}. In particular, they do not change under
coordinate change or when topologically irrelevant pairs are omitted.

In this splice diagram the arrowhead represents the component of the
link and the weights along and adjacent to the path to this arrowhead
give the sequence of cabling pairs.

To understand the placement issue for more than one branch, suppose,
for example, that $k=3$ above, so there are just three characteristic
pairs. Suppose also that our curve has a second branch whose first
pair is also $(p_1,s_1)$ but whose later pairs are different, say
$(p_2',s_2'), (p_3',s_3')$.  Assume also that $\frac{p_2'}{s_2'}$ is
the larger of $\frac{p_2}{s_2}$ and $\frac{p_2'}{s_2'}$. The splice
diagram may then be
$$
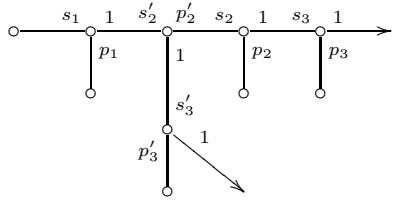
\xymatrix@R=9pt@C=24pt@M=0pt@W=0pt@H=0pt{
\circ\lineto[r]^(.75){s_1}&
\circ\lineto[dd]^(.35){p_1}\lineto[r]^(.25){1}^(.75){s_2'}&
\circ\lineto[ddd]^(.25){1}^(.75){s_3'}\lineto[r]^(.25){p_2'}^(.75){s_2}&
\circ\lineto[dd]^(.35){p_2}\lineto[r]^(.25){1}^(.75){s_3}&
\circ\lineto[dd]^(.35){p_3}\ar[r]^(.25){1}&\\ \\ 
&\circ&&\circ&\circ\\&&\circ\lineto[dd]_(.35){p_3'}\ar[ddr]^(.35){1}
\\ \\
&&\circ&&}
$$\iffalse
or $$\xymatrix@R=9pt@C=24pt@M=0pt@W=0pt@H=0pt{
\circ\lineto[r]^(.75){s_1}&
\circ\lineto[dd]^(.35){p_1}\lineto[r]^(.25){1}^(.75){s_2}&
\circ\lineto[ddd]^(.25){1}^(.75){s_3}\lineto[r]^(.25){p_2}^(.75){s_2'}&
\circ\lineto[dd]^(.35){p_2'}\lineto[r]^(.25){1}^(.75){s_3'}&
\circ\lineto[dd]^(.35){p_3'}\ar[r]^(.25){1}&\\ \\ 
&\circ&&\circ&\circ\\&&\circ\lineto[dd]_(.35){p_3}\ar[ddr]^(.35){1}
\\ \\
&&\circ&&}
$$\fi
or 
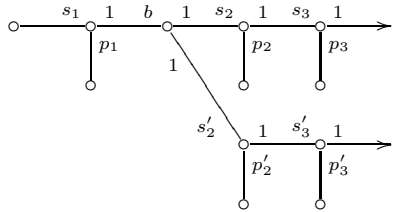
$$\xymatrix@R=18pt@C=24pt@M=0pt@W=0pt@H=0pt{
\circ\lineto[r]^(.75){s_1}&
\circ\lineto[d]^(.35){p_1}\lineto[r]^(.25){1}^(.75){b}&
\circ\lineto[ddr]_(.25){1}_(.75){s_2'}\lineto[r]^(.25){1}^(.75){s_2}&
\circ\lineto[d]^(.35){p_2}\lineto[r]^(.25){1}^(.75){s_3}&
\circ\lineto[d]^(.35){p_3}\ar[r]^(.25){1}&\\
&\circ&&\circ&\circ&\\%\hbox to 0pt{\qquad$(b>s_1p_1)$\hss}\\
&&&\circ\lineto[d]^(.35){p_2'}\lineto[r]^(.25){1}^(.75){s_3'}&
\circ\lineto[d]^(.35){p_3'}\ar[r]^(.25){1}&\\
&&&\circ&\circ&
}
$$or
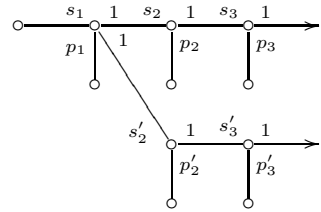
$$\xymatrix@R=18pt@C=24pt@M=0pt@W=0pt@H=0pt{
\circ\lineto[r]^(.75){s_1}&
\circ\lineto[d]_(.4){p_1}\lineto[r]^(.25){1}^(.75){s_2}
\lineto[ddr]^(.2){1}_(.8){s_2'}&
\circ\lineto[d]^(.35){p_2}\lineto[r]^(.25){1}^(.75){s_3}&
\circ\lineto[d]^(.35){p_3}\ar[r]^(.25){1}&\\
&\circ&\circ&\circ\\
&&
\circ\lineto[d]^(.35){p_2'}\lineto[r]^(.25){1}^(.75){s_3'}&
\circ\lineto[d]^(.35){p_3'}\ar[r]^(.25){1}&\\
&&\circ&\circ}
$$
or 
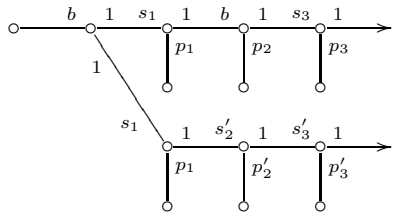
$$\xymatrix@R=18pt@C=24pt@M=0pt@W=0pt@H=0pt{
\circ\lineto[r]^(.75){b}&
\circ\lineto[ddr]_(.25){1}_(.75){s_1}\lineto[r]^(.25){1}^(.75){s_1}&
\circ\lineto[d]^(.35){p_1}\lineto[r]^(.25){1}^(.75){b}&
\circ\lineto[d]^(.35){p_2}\lineto[r]^(.25){1}^(.75){s_3}&
\circ\lineto[d]^(.35){p_3}\ar[r]^(.25){1}&\\
&&\circ&\circ&\circ\\
&&\circ\lineto[d]^(.35){p_1}\lineto[r]^(.25){1}^(.75){s_2'}
&\circ\lineto[d]^(.35){p_2'}\lineto[r]^(.25){1}^(.75){s_3'}&
\circ\lineto[d]^(.35){p_3'}\ar[r]^(.25){1}&\\
&&\circ&\circ&\circ&
}
$$or
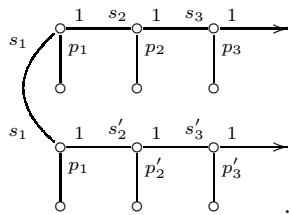
$$\xymatrix@R=18pt@C=24pt@M=0pt@W=0pt@H=0pt{
&
%\circ\lineto[ddr]_(.25){1}_(.75){s_1}\lineto[r]^(.25){1}^(.75){s_1}
&\circ\ar@/_14pt/@{-}[dd]_(.2){s_1}_(.8){s_1}
\lineto[d]^(.35){p_1}\lineto[r]^(.25){1}^(.75){s_2}&
\circ\lineto[d]^(.35){p_2}\lineto[r]^(.25){1}^(.75){s_3}&
\circ\lineto[d]^(.35){p_3}\ar[r]^(.25){1}&\\
&&\circ&\circ&\circ\\
&&\circ\lineto[d]^(.35){p_1}\lineto[r]^(.25){1}^(.75){s_2'}
&\circ\lineto[d]^(.35){p_2'}\lineto[r]^(.25){1}^(.75){s_3'}&
\circ\lineto[d]^(.35){p_3'}\ar[r]^(.25){1}&\\
&&\circ&\circ&\circ&.
}
$$\\

There are conditions on 
%the edge weights of a 
splice diagrams that are necessary and sufficient for a diagram to
occur for the link of some plane curve singularity (and be the unique
smallest diagram representing its topology). The edge weights
%in a general splice diagram are weights on edges that 
%occur adjacent to nodes (vertices of valence $>1$), and the weights
around any %particular 
node (vertex of valence $>1$) are pairwise coprime.  Modification
procedures on a splice diagram can sometimes create an edge weight
adjacent to a leaf (vertex of valence 1), but if this happens the
weight should just be deleted.  The other conditions are:
%for splice diagrams of plane curve singularities are:
\begin{itemize}
\item all edge weights are positive;
\item each \emph{edge determinant} is positive (the edge determinant,
  defined for each edge connecting two nodes, is the product of the
  two weights on the edge minus the product of the weights directly
  adjacent to the edge);
\item An edge to a leaf should not have weight 1 (if it does, remove
  the edge);
\item no vertex should have valence 2 (eliminate such a vertex by
  replacing it and its two adjacent edges by a single edge);
\item if all arrowheads are replaced by vertices, the diagram should
  collapse to a single vertex or single edge using the moves just
  described.
\end{itemize}
The linking number of two link components is particularly easy to
compute in terms of the splice diagram: it is the product of the
weights adjacent to but not on the path that connects the
corresponding arrowheads. For example, in the first of the above
diagrams it is $s_2'p_2p_3p_3'$. The positive edge determinant
condition implies immediately that this linking number is
strictly decreasing as one runs through the five possible placements
of the two branches in the above example.  It is easy to turn this
into a quick general splice diagram proof of Theorem \ref{th:1}.

With different conditions, splice diagrams can encode many other
objects of interest.  For example, the last condition above just
assures that we are looking at a point of a curve at a non-singular
point of a surface.
One advantage of splice diagrams is that
invariants such as fundamental group, Alexander polynomial
(single-variable and multi-variable), Milnor fiber, value semigroup,
etc., can be computed quite easily and uniformly from the splice
diagram in any situation where the invariant makes sense. Such
situations include the study of the global topology of plane curves
(work of Neumann, Neumann and Norbury, Pierrette Cassou-Nogu\`es, and
others; here the r\^ole of Milnor fiber is played by the generic fiber
of the defining polynomial of the curve), and the study of surface
singularities with homology sphere links. Recently, splice diagrams in
which the coprimality condition is relaxed have been used in the study
of universal abelian covers of surface singularities (Neumann
and Wahl).

\smallskip We have surveyed here the \emph{topology} of plane curve
singularities and intentionally not ventured into the large and active
literature on algebraic/analytic aspects such as deformation and
moduli spaces, curves over fields of finite characteristic, etc. Even
with this restriction we have had to leave much out. The 1981 book
\cite{brieskorn} of Brieskorn and Kn\"orrer is a delightful and
readable survey from ancient times to 1980. But the subject has not
stopped there. Very recent papers include A'Campo's beautiful
construction of the Milnor fibration from a real morsification
\cite{a'campo}, and an intriguing and surprising geometric formula for
the Alexander polynomial in \cite{zade}.

\section{Topology of hypersurface singularities in higher dimension}

The study of the topology of isolated singularities of complex
hypersurfaces in higher dimensions received a considerable boost in
the late 1960's from Brieskorn's construction of exotic spheres as
singularity links (of what are now known as Brieskorn-Pham
singularities), and from Milnor's monograph \cite{milnor}. Durfee
\cite{durfee99} gives an excellent history of this period. 
Milnor's fibration theorem is now a fundamental tool in
the subject. It says that the link of an isolated hypersurface
singularity is a fibered link (that is, the complement of the link can
be fibered over $S^1$ with fibers which are the interiors of
submanifolds of the sphere with the link as boundary). Milnor proved
that on the standard sphere $\partial D^{2n}(\epsilon)$ the fibration
of the complement of the link is given by $f/|f|$, but this is rarely
needed, so we will sketch a version of Milnor's proof that omits this
fact, but has some of the spirit of K\"ahler's paper.

Suppose $f:\C^n\to \C$ is such that $f(0)=0$ and $f$ has an isolated
singularity at $0\in\C^n$. Then for $\delta$ and $\epsilon$
sufficiently small, and $\delta <<\epsilon$, a vector field argument
shows that $D:=f^{-1}(D^2(\delta))\cap D^{2n}(\epsilon)$ is
isotopically equivalent to the ball $D^{2n}(\epsilon)$. We can thus
consider the link of the singularity in the boundary of this (somewhat
twisted) ``rectangular ball'' $D$. The function $f$ restricted to
$\partial D$ makes the desired fibered structure evident.

The fiber of Milnor's fibration is highly connected, with homology
only in its middle dimension. It is a Seifert surface for the link,
and the Seifert linking form with respect to it is a natural algebraic
invariant. (The Seifert form evaluates linking numbers of cycles in
the Seifert surface with cycles in a parallel copy of the surface; it
is a non-singular integral bilinear form.) Durfee observed in
\cite{durfee} that work of Levine \cite{levine} implies that the
Seifert form is a complete invariant for the topology of a link of an
isolated hypersurface in $\C^n$ if $n>3$.

In \cite{eisenbud-neumann} it was asked if the Seifert form is a
complete topological invariant for plane curve singularities ($n=2$).
Counterexamples were found by Du Bois and Michel \cite{dubois}, and
used by Artal-Bartolo \cite{artal} to give a counterexample also for
$n=3$.

Isolated hypersurface singularities in $\C^3$ thus remain the
topologically least well understood.  The link is a 3-manifold in
$S^5$.  Although one knows what 3-manifolds can be links of isolated
surface singularities (by the work of Grauert combined with standard
3-manifold theory---see, e.g., \cite{neumann81}), it is not known
which of them occur for {hypersurface} singularities, and the possible
embeddings in $S^5$ as links of hypersurface singularities are even
less understood.  If the $3$-manifold is $S^3$ then there is no
singularity and the embedding is standard, but other $3$-manifolds may
have several embeddings as singularity links.

This is not to say that the topology \emph{is} understood in
dimensions $n>3$.  The Durfee-Levine theorem says that the Seifert
form tells all, but it is unknown what forms $L$ actually occur as
Seifert forms of singularity links in dimension $n$. Some restrictions
are known. For example, there is a basis that makes the form $L$ upper
triangular with diagonal entries $(-1)^{\frac{n(n-1)}2}$ (Durfee
\cite{durfee}).  The eigenvalues of $L^tL^{-1}$ are roots of unity
with maximal Jordan block size $n$ by the monodromy theorem of
Grothendieck and Deligne (see \cite{vandoorn} for a slight
sharpening).  It is conjectured (problem 3.31 of %Kirby's Problem List
\cite{kirby}, attributed to Durfee \cite{durfee78}) that
$(-1)^{\frac{n(n-1)}2}(L+L^t)$ always has positive signature. Durfee's
original conjecture was only for $n=2$, but even this is unknown.

Since ``suspending'' a singularity (replacing the hypersurface
$f(x_1,\dots,x_n)=0$ by $f(x_1,\dots,x_n)+x_{n+1}^2 =0$) just
multiplies the Seifert form by $(-1)^n$, if we adjust sign of the
Seifert form by replacing $L$ by $(-1)^{\frac{n(n-1)}2}L$, then
the set of realized forms grows with dimension. This %interesting
graded set of forms is closed, in a graded sense, under tensor
product.  It fully describes the topology of all isolated hypersurface
singularities in ambient dimensions $>3$, but it remains mysterious.

\iffalse The tensor product of forms realized at levels $r$ and $s$ is
realised at level $r+s$ by the Seifert form version of the
Thom-Sebastiani theorem (proved independently by Deligne, Gabrielov,
Sakamoto).  \fi

\end{document}